\documentclass[a4paper,10pt]{article}
\usepackage{amsmath,amsthm,amssymb,amscd,layout}   
\usepackage{bm,color,enumerate} 
\usepackage{amsfonts}
\usepackage{mathrsfs}
\usepackage[dvips]{graphicx}

\newtheorem{theorem}{Theorem}[section]
\newtheorem{lemma}[theorem]{Lemma}
\newtheorem{proposition}[theorem]{Proposition}

\theoremstyle{definition}
\newtheorem{definition}[theorem]{Definition}
\newtheorem{example}[theorem]{Example}

\theoremstyle{remark}
\newtheorem{remark}[theorem]{Remark}

\theoremstyle{conjecture}

\numberwithin{equation}{section}
\newcommand{\hilsp}{\mathscr{H}}

\newcommand{\sfrac}[2]{
\mbox{\small $\dfrac{#1}{#2}$}
}

\newcommand{\crzer}[3]{
\setlength{\unitlength}{1.6pt}
\begin{picture}(20,55)(0,-8)
   \thinlines
   \multiput(  0, 0)(20,0){2}{\line(0,1){40}}
   \multiput(  0, 0)(0,40){2}{\line(1,0){20}}
   \put(  0, 40){\makebox(20, 7){$ #1 $}}
   \put(  0, -8){\makebox(20, 8){$ #2 $}}
   \put(  0,-14){\makebox(20, 6){$\scriptstyle{#3}$}}
   \thicklines
   \put( 19, 1){\oval(18,18)[tl]}
   \put( 19,10){\line(1,0){1}}
   \put( 10, 1){\circle*{2}}
\end{picture}
}
\newcommand{\crone}[3]{
\setlength{\unitlength}{1.6pt}
\begin{picture}(20,55)(0,-8)
   \thinlines
   \multiput(  0, 0)(20,0){2}{\line(0,1){40}}
   \multiput(  0, 0)(0,40){2}{\line(1,0){20}}
   \put(  0, 40){\makebox(20, 7){$  #1 $}}
   \put(  0, -8){\makebox(20, 8){$  #2 $}}
   \put(  0,-14){\makebox(20, 6){$\scriptstyle{#3}$}}
   \thicklines
   \put( 19, 1){\oval(18,18)[tl]}
   \put( 19,10){\line(1,0){1}}
   \put( 10, 1){\circle*{2}}
   \qbezier( 0,10)( 4,10)(10,14)
   \qbezier(19,18)(16,18)(10,14)
   \put(19,18){\line(1,0){1}}
\end{picture}
}
\newcommand{\crtwo}[3]{
\setlength{\unitlength}{1.6pt}
\begin{picture}(20,55)(0,-8)
   \thinlines
   \multiput(  0, 0)(20,0){2}{\line(0,1){40}}
   \multiput(  0, 0)(0,40){2}{\line(1,0){20}}
   \put(  0, 40){\makebox(20, 7){$  #1 $}}
   \put(  0, -8){\makebox(20, 8){$  #2 $}}
   \put(  0,-14){\makebox(20, 6){$\scriptstyle{#3}$}}
   \thicklines
   \put( 19, 1){\oval(18,18)[tl]}
   \put( 19,10){\line(1,0){1}}
   \put( 10, 1){\circle*{2}}
   \qbezier( 0,10)( 4,10)(10,14)
   \qbezier(19,18)(16,18)(10,14)
   \put(19,18){\line(1,0){1}}
   \qbezier( 0,18)( 4,18)(10,22)
   \qbezier(19,26)(16,26)(10,22)
   \put(19,26){\line(1,0){1}}
\end{picture}
}
\newcommand{\crnth}[3]{
\setlength{\unitlength}{1.6pt}
\begin{picture}(20,55)(0,-8)
   \thinlines
   \multiput(  0, 0)(20,0){2}{\line(0,1){40}}
   \multiput(  0, 0)(0,40){2}{\line(1,0){20}}
   \put(  0, 40){\makebox(20, 7){$ #1 $}}
   \put(  0, -8){\makebox(20, 8){$ #2 $}}
   \put(  0,-14){\makebox(20, 6){$\scriptstyle{#3}$}}
   \thicklines
   \put( 19, 1){\oval(18,18)[tl]}
   \put( 19,10){\line(1,0){1}}
   \put( 10, 1){\circle*{2}}
   \qbezier( 0,10)( 4,10)(10,14)
   \qbezier(19,18)(16,18)(10,14)
   \put(19,18){\line(1,0){1}}
   \put(10,24){\makebox(0,0){$ \vdots $}}
   \qbezier( 0,26)( 4,26)(10,30)
   \qbezier(19,34)(16,34)(10,30)
   \put(19,34){\line(1,0){1}}
\end{picture}
}
\newcommand{\anone}[3]{
\setlength{\unitlength}{1.6pt}
\begin{picture}(20,55)(0,-8)
   \thinlines
   \multiput(  0, 0)(20,0){2}{\line(0,1){40}}
   \multiput(  0, 0)(0,40){2}{\line(1,0){20}}
   \put(  0, 40){\makebox(20, 7){$ #1 $}}
   \put(  0, -8){\makebox(20, 8){$ #2 $}}
   \put(  0,-14){\makebox(20, 6){$\scriptstyle{#3}$}}
   \thicklines
   \put(  1, 1){\oval(18,18)[tr]}
   \put(  0,10){\line(1,0){1}}
   \put( 10, 1){\circle*{2}}
\end{picture}
}
\newcommand{\antwo}[3]{
\setlength{\unitlength}{1.6pt}
\begin{picture}(20,55)(0,-8)
   \thinlines
   \multiput(  0, 0)(20,0){2}{\line(0,1){40}}
   \multiput(  0, 0)(0,40){2}{\line(1,0){20}}
   \put(  0, 40){\makebox(20, 7){$ #1 $}}
   \put(  0, -8){\makebox(20, 8){$ #2 $}}
   \put(  0,-14){\makebox(20, 6){$\scriptstyle{#3}$}}
   \thicklines
   \put(  1, 1){\oval(18,18)[tr]}
   \put(  0,10){\line(1,0){1}}
   \put( 10, 1){\circle*{2}}
   \qbezier( 0,18)( 4,18)(10,14)
   \qbezier(19,10)(16,10)(10,14)
   \put(19,10){\line(1,0){1}}
\end{picture}
}
\newcommand{\anthr}[3]{
\setlength{\unitlength}{1.6pt}
\begin{picture}(20,55)(0,-8)
   \thinlines
   \multiput(  0, 0)(20,0){2}{\line(0,1){40}}
   \multiput(  0, 0)(0,40){2}{\line(1,0){20}}
   \put(  0, 40){\makebox(20, 7){$ #1 $}}
   \put(  0, -8){\makebox(20, 8){$ #2 $}}
   \put(  0,-14){\makebox(20, 6){$\scriptstyle{#3}$}}
   \thicklines
   \put(  1, 1){\oval(18,18)[tr]}
   \put(  0,10){\line(1,0){1}}
   \put( 10, 1){\circle*{2}}
   \qbezier( 0,18)( 4,18)(10,14)
   \qbezier(19,10)(16,10)(10,14)
   \put(19,10){\line(1,0){1}}
   \qbezier( 0,26)( 4,26)(10,22)
   \qbezier(19,18)(16,18)(10,22)
   \put(19,18){\line(1,0){1}}
\end{picture}
}
\newcommand{\annth}[3]{
\setlength{\unitlength}{1.6pt}
\begin{picture}(20,55)(0,-8)
   \thinlines
   \multiput(  0, 0)(20,0){2}{\line(0,1){40}}
   \multiput(  0, 0)(0,40){2}{\line(1,0){20}}
   \put(  0, 40){\makebox(20, 7){$ #1 $}}
   \put(  0, -8){\makebox(20, 8){$ #2 $}}
   \put(  0,-14){\makebox(20, 6){$\scriptstyle{#3}$}}
   \thicklines
   \put(  1, 1){\oval(18,18)[tr]}
   \put(  0,10){\line(1,0){1}}
   \put( 10, 1){\circle*{2}}
   \qbezier( 0,18)( 4,18)(10,14)
   \qbezier(19,10)(16,10)(10,14)
   \put(19,10){\line(1,0){1}}
   \put(10,24){\makebox(0,0){$ \vdots $}}
   \qbezier( 0,34)( 4,34)(10,30)
   \qbezier(19,26)(16,26)(10,30)
   \put(19,26){\line(1,0){1}}
\end{picture}
}
\newcommand{\nuone}[3]{
\setlength{\unitlength}{1.6pt}
\begin{picture}(20,55)(0,-8)
   \thinlines
   \multiput(  0, 0)(20,0){2}{\line(0,1){40}}
   \multiput(  0, 0)(0,40){2}{\line(1,0){20}}
   \put(  0, 40){\makebox(20, 7){$ #1 $}}
   \put(  0, -8){\makebox(20, 8){$ #2 $}}
   \put(  0,-14){\makebox(20, 6){$\scriptstyle{#3}$}}
   \thicklines
   \put(  1, 1){\oval(18,18)[tr]}
   \put(  0,10){\line(1,0){1}}
   \put( 10, 1){\circle*{2}}
   \put( 19, 1){\oval(18,18)[tl]}
   \put( 19,10){\line(1,0){1}}
   \put( 10, 1){\circle*{2}}
\end{picture}
}
\newcommand{\nutwo}[3]{
\setlength{\unitlength}{1.6pt}
\begin{picture}(20,55)(0,-8)
   \thinlines
   \multiput(  0, 0)(20,0){2}{\line(0,1){40}}
   \multiput(  0, 0)(0,40){2}{\line(1,0){20}}
   \put(  0, 40){\makebox(20, 7){$ #1 $}}
   \put(  0, -8){\makebox(20, 8){$ #2 $}}
   \put(  0,-14){\makebox(20, 6){$\scriptstyle{#3}$}}
   \thicklines
   \put( 19, 1){\oval(18,18)[tl]}
   \put( 19,10){\line(1,0){1}}
   \put( 10, 1){\circle*{2}}
   \put(  1, 1){\oval(18,18)[tr]}
   \put(  0,10){\line(1,0){1}}
   \put( 10, 1){\circle*{2}}
   \put( 0, 18){\line(1,0){20}}
\end{picture}
}
\newcommand{\nuthr}[3]{
\setlength{\unitlength}{1.6pt}
\begin{picture}(20,55)(0,-8)
   \thinlines
   \multiput(  0, 0)(20,0){2}{\line(0,1){40}}
   \multiput(  0, 0)(0,40){2}{\line(1,0){20}}
   \put(  0, 40){\makebox(20, 7){$ #1 $}}
   \put(  0, -8){\makebox(20, 8){$ #2 $}}
   \put(  0,-14){\makebox(20, 6){$\scriptstyle{#3}$}}
   \thicklines
   \put( 19, 1){\oval(18,18)[tl]}
   \put( 19,10){\line(1,0){1}}
   \put( 10, 1){\circle*{2}}
   \put(  1, 1){\oval(18,18)[tr]}
   \put(  0,10){\line(1,0){1}}
   \put( 10, 1){\circle*{2}}
   \put( 0, 18){\line(1,0){20}}
   \put( 0, 26){\line(1,0){20}}
\end{picture}
}
\newcommand{\nunth}[3]{
\setlength{\unitlength}{1.6pt}
\begin{picture}(20,55)(0,-8)
   \thinlines
   \multiput(  0, 0)(20,0){2}{\line(0,1){40}}
   \multiput(  0, 0)(0,40){2}{\line(1,0){20}}
   \put(  0, 40){\makebox(20, 7){$ #1 $}}
   \put(  0, -8){\makebox(20, 8){$ #2 $}}
   \put(  0,-14){\makebox(20, 6){$\scriptstyle{#3}$}}
   \thicklines
   \put( 19, 1){\oval(18,18)[tl]}
   \put( 19,10){\line(1,0){1}}
   \put( 10, 1){\circle*{2}}
   \put(  1, 1){\oval(18,18)[tr]}
   \put(  0,10){\line(1,0){1}}
   \put( 10, 1){\circle*{2}}
   \put( 0, 18){\line(1,0){20}}
   \put(10,27){\makebox(0,0){$ \vdots $}}
   \put( 0, 34){\line(1,0){20}}
\end{picture}
}
\newcommand{\kozer}[3]{
\setlength{\unitlength}{1.6pt}
\begin{picture}(20,55)(0,-8)
   \thinlines
   \multiput(  0, 0)(20,0){2}{\line(0,1){40}}
   \multiput(  0, 0)(0,40){2}{\line(1,0){20}}
   \put(  0, 40){\makebox(20, 7){$ #1 $}}
   \put(  0, -8){\makebox(20, 8){$ #2 $}}
   \put(  0,-14){\makebox(20, 6){$\scriptstyle{#3}$}}
   \thicklines
   \put( 10, 0){\line(0,1){6}}
   \put( 10, 1){\circle*{2}}
\end{picture}
}
\newcommand{\koone}[3]{
\setlength{\unitlength}{1.6pt}
\begin{picture}(20,55)(0,-8)
   \thinlines
   \multiput(  0, 0)(20,0){2}{\line(0,1){40}}
   \multiput(  0, 0)(0,40){2}{\line(1,0){20}}
   \put(  0, 40){\makebox(20, 7){$ #1 $}}
   \put(  0, -8){\makebox(20, 8){$ #2 $}}
   \put(  0,-14){\makebox(20, 6){$\scriptstyle{#3}$}}
   \thicklines
   \put( 10, 0){\line(0,1){6}}
   \put( 10, 1){\circle*{2}}
   \put(  0, 10){\line(1,0){20}}
\end{picture}
}
\newcommand{\kotwo}[3]{
\setlength{\unitlength}{1.6pt}
\begin{picture}(20,55)(0,-8)
   \thinlines
   \multiput(  0, 0)(20,0){2}{\line(0,1){40}}
   \multiput(  0, 0)(0,40){2}{\line(1,0){20}}
   \put(  0, 40){\makebox(20, 7){$ #1 $}}
   \put(  0, -8){\makebox(20, 8){$ #2 $}}
   \put(  0,-14){\makebox(20, 6){$\scriptstyle{#3}$}}
   \thicklines
   \put( 10, 0){\line(0,1){6}}
   \put( 10, 1){\circle*{2}}
   \put(  0, 10){\line(1,0){20}}
   \put(  0, 18){\line(1,0){20}}
\end{picture}
}
\newcommand{\konth}[3]{
\setlength{\unitlength}{1.6pt}
\begin{picture}(20,55)(0,-8)
   \thinlines
   \multiput(  0, 0)(20,0){2}{\line(0,1){40}}
   \multiput(  0, 0)(0,40){2}{\line(1,0){20}}
   \put(  0, 40){\makebox(20, 7){$ #1 $}}
   \put(  0, -8){\makebox(20, 8){$ #2 $}}
   \put(  0,-14){\makebox(20, 6){$\scriptstyle{#3}$}}
   \thicklines
   \put( 10, 0){\line(0,1){6}}
   \put( 10, 1){\circle*{2}}
   \put(  0, 10){\line(1,0){20}}
   \put(10,19){\makebox(0,0){$ \vdots $}}
   \put(  0, 26){\line(1,0){20}}
\end{picture}
}
\newcommand{\rsideihi}[1]{
\setlength{\unitlength}{1.6pt}
\begin{picture}(10,55)(0,-8)
   \put(0,25){{$\left. \makebox(0,10){} \right\} \scriptstyle{#1} $ }}
\end{picture}
}
\newcommand{\rsideilo}[1]{
\setlength{\unitlength}{1.6pt}
\begin{picture}(10,55)(0,-8)
   \put(0,17){{$\left. \makebox(0,10){} \right\} \scriptstyle{#1} $ }}
\end{picture}
}
\newcommand{\lsideihi}[1]{
\setlength{\unitlength}{1.6pt}
\begin{picture}(7,55)(0,-8)
   \put(0,25){{$ \scriptstyle{#1} \left\{ \makebox(0,11){} \right. $ }}
\end{picture}
}
\newcommand{\drawarc}[2]{
\ifnum#1=1
  \ifnum#2=1
      \put(  0, 0){\line(0,1){5}}
  \fi
  \ifnum#2=2
      \qbezier(0, 0)( 4, 8)( 8, 0)
  \fi
  \ifnum#2=3
      \qbezier(0, 0)( 8,16)(16, 0)
  \fi
  \ifnum#2=4
      \qbezier(0, 0)(12,24)(24, 0)
  \fi
  \ifnum#2=5
      \qbezier(0, 0)(16,32)(32, 0)
  \fi
  \ifnum#2=6
      \qbezier(0, 0)(20,40)(40, 0)
  \fi
  \ifnum#2=7
      \qbezier(0, 0)(24,48)(48, 0)
  \fi
  \ifnum#2=8
      \qbezier(0, 0)(28,56)(56, 0)
  \fi
  \ifnum#2=9
      \qbezier(0, 0)(32,64)(64, 0)
  \fi
  \ifnum#2=10
      \qbezier(0, 0)(36,72)(72, 0)
  \fi
\fi
\ifnum#1=2
  \ifnum#2=2
      \put(  8, 0){\line(0,1){5}}
  \fi
  \ifnum#2=3
      \qbezier( 8, 0)(12, 8)(16, 0)
  \fi
  \ifnum#2=4
      \qbezier( 8, 0)(16,16)(24, 0)
  \fi
  \ifnum#2=5
      \qbezier( 8, 0)(20,24)(32, 0)
  \fi
  \ifnum#2=6
      \qbezier( 8, 0)(24,32)(40, 0)
  \fi
  \ifnum#2=7
      \qbezier( 8, 0)(28,40)(48, 0)
  \fi
  \ifnum#2=8
      \qbezier( 8, 0)(32,48)(56, 0)
  \fi
  \ifnum#2=9
      \qbezier( 8, 0)(36,64)(64, 0)
  \fi
  \ifnum#2=10
      \qbezier( 8, 0)(40,72)(72, 0)
  \fi
\fi
\ifnum#1=3
  \ifnum#2=3
      \put( 16, 0){\line(0,1){5}}
  \fi
  \ifnum#2=4
      \qbezier(16, 0)(20, 8)(24, 0)
  \fi
  \ifnum#2=5
      \qbezier(16, 0)(24,16)(32, 0)
  \fi
  \ifnum#2=6
      \qbezier(16, 0)(28,24)(40, 0)
  \fi
  \ifnum#2=7
      \qbezier(16, 0)(32,32)(48, 0)
  \fi
  \ifnum#2=8
      \qbezier(16, 0)(36,40)(56, 0)
  \fi
  \ifnum#2=9
      \qbezier(16, 0)(40,64)(64, 0)
  \fi
  \ifnum#2=10
      \qbezier(16, 0)(44,72)(72, 0)
  \fi
\fi
\ifnum#1=4
  \ifnum#2=4
      \put( 24, 0){\line(0,1){5}}
  \fi
  \ifnum#2=5
      \qbezier(24, 0)(28, 8)(32, 0)
  \fi
  \ifnum#2=6
      \qbezier(24, 0)(32,16)(40, 0)
  \fi
  \ifnum#2=7
      \qbezier(24, 0)(36,24)(48, 0)
  \fi
  \ifnum#2=8
      \qbezier(24, 0)(40,32)(56, 0)
  \fi
  \ifnum#2=9
      \qbezier(24, 0)(44,40)(64, 0)
  \fi
  \ifnum#2=10
      \qbezier(24, 0)(48,48)(72, 0)
  \fi
\fi
\ifnum#1=5
  \ifnum#2=5
      \put( 32, 0){\line(0,1){5}}
  \fi
  \ifnum#2=6
      \qbezier(32, 0)(36, 8)(40, 0)
  \fi
  \ifnum#2=7
      \qbezier(32, 0)(40,16)(48, 0)
  \fi
  \ifnum#2=8
      \qbezier(32, 0)(44,24)(56, 0)
  \fi
  \ifnum#2=9
      \qbezier(32, 0)(48,32)(64, 0)
  \fi
  \ifnum#2=10
      \qbezier(32, 0)(52,40)(72, 0)
  \fi
\fi
\ifnum#1=6
  \ifnum#2=6
      \put( 40, 0){\line(0,1){5}}
  \fi
  \ifnum#2=7
      \qbezier(40, 0)(44, 8)(48, 0)
  \fi
  \ifnum#2=8
      \qbezier(40, 0)(48,16)(56, 0)
  \fi
  \ifnum#2=9
      \qbezier(40, 0)(52,24)(64, 0)
  \fi
  \ifnum#2=10
      \qbezier(40, 0)(56,32)(72, 0)
  \fi
\fi
\ifnum#1=7
  \ifnum#2=7
      \put( 48, 0){\line(0,1){5}}
  \fi
  \ifnum#2=8
      \qbezier(48, 0)(52, 8)(56, 0)
  \fi
  \ifnum#2=9
      \qbezier(48, 0)(56,40)(64, 0)
  \fi
  \ifnum#2=10
      \qbezier(48, 0)(60,48)(72, 0)
  \fi
\fi
\ifnum#1=8
  \ifnum#2=8
      \put( 56, 0){\line(0,1){5}}
  \fi
  \ifnum#2=9
      \qbezier(56, 0)(60,40)(64, 0)
  \fi
  \ifnum#2=10
      \qbezier(56, 0)(64,48)(72, 0)
  \fi
\fi
\ifnum#1=9
  \ifnum#2=9
      \put( 64, 0){\line(0,1){5}}
  \fi
  \ifnum#2=10
      \qbezier(64, 0)(68,40)(64, 0)
  \fi
\fi
\ifnum#1=10
  \ifnum#2=10
      \put( 72, 0){\line(0,1){5}}
  \fi
\fi
}

\title{\bf Combinatorial Aspects of\\
Weighted Free Poisson Random 
Variables\footnote{Submitted: 
August 31, 2023.  Accepted: January 5, 2024}}
\author{Nobuhiro ASAI\footnote{Supported by 
JSPS KAKENHI Grant Numbers JP20K03652. }\\
Department of Mathematics, \\
Aichi University of Education, \\ 
Hirosawa 1, Igaya,  \\
Kariya 448-8542, Japan. \\
E-mail:  {\tt nasai[at]auecc.aichi-edu.ac.jp} \\ 
and \\ 
Hiroaki YOSHIDA\footnote{Supported by 
JSPS KAKENHI Grant Numbers JP20K03649.}\\
Department of Information Sciences, \\
Ochanomizu University, \\
Tokyo 112-8610, Japan. \\
E-mail: {\tt yoshida[at]is.ocha.ac.jp}
}

\date{}

\begin{document}

\maketitle


\begin{abstract}
This paper will be devoted to study 
weighted (deformed) free Poisson random variables
from the viewpoint of orthogonal polynomials and statistics of non-crossing partitions. 
A family of weighted (deformed) free Poisson random variables will 
be defined in a sense by the sum of  weighted (deformed) free creation, annihilation, 
scalar, and intermediate operators with certain parameters 
on a weighted (deformed) free Fock space together with the vacuum expectation.
We shall provide a combinatorial moment formula of non-commutative 
Poisson random variables.
This formula gives us a very nice combinatorial interpretation to 
two parameters of weights.  One can see that 
the deformation treated in this paper interpolates
free and boolean Poisson random variables,  their distributions and moments,
and yields some conditionally free Poisson distribution by taking limit of 
the parameter.


\smallskip 


\end{abstract}

%
\section{Introduction}\label{sec1}
From the probabilistic point of view,  
the Poisson distribution in addition to Gaussian is one of fundamental
objects to be considered because arbitrary infinitely divisible distributions
can be constructed by using the Gaussian and Poisson distributions due to 
L\`evy-Khintchine representation.  

From the non-commutative probabilistic point of view, 
based on the conditionally free product of states \cite{BLS96},  
a large class of deformed free convolution, so-called $\varDelta$-convolution,
was introduced.   In \cite{B01}\cite{Y02},  
the $s$-free convolution was treated as an interesting 
example of the $\varDelta$-convolution.
The $s$-deformation is an attempt to realize an interpolation between the free product 
of states and boolean (regular free) product of the states.   
See \cite{B01}\cite{VDN92} and references therein.
In \cite{Y03},  the $\varDelta$-deformed moment-cumulant formula  
and $\varDelta$-free Gaussian and Poisson distributions
were obtained associated with a certain very general weight function 
for set partition statistics on non-crossing pair partitions.  
One can imagine easily that
the sum of $s$-creation and annihilation operators 
acting on the $s$-free Fock 
space (a certain weighted full Fock space)
plays a role of the $s$-free analogue of Gaussian field operator and 
can be seen as a realization of non-commutative random variable associated with the 
$s$-free Gaussian distribution.   See \cite{BKW06} for the $r$-free Gaussian case.

The main purpose of this paper is to 
examine the $s$-free Poisson counterpart 
by combinatorial consideration.   
We note that our partition statistics are based on 
the last and intermediate elements of the block on the non-crossing partitions
and hence are more natural than those in \cite{Y02}\cite{Y03}
to consider the $s$-free Poisson part on the $s$-free Fock space.

This paper will be organized as follows.
Firstly,  we shall introduce an analogue of the Poisson type random variable, namely,   
the $s$-weighted (deformed) free Poisson random variable as a non-commutative random 
variable accompanied with the weighted (deformed) free Poisson distribution in Section \ref{sec3}. 
Although our weighted (deformed) free Poisson random variable  
has the similar form as of the free case in \cite{Sp90},  
counterparts of gauge and identity operators 
should be replaced by the intermediate and scalar operators given in Section \ref{sec2}, respectively. 
The reason of these replacements comes from 
combinatorial consideration on non-crossing pair partitions
discussed in Section \ref{sec4} and Section \ref{sec5}. 
Secondly, we shall give the recurrence formula for the orthogonal polynomials with respect to 
the $s$-free Poisson distribution in Section \ref{sec3}. 
In order to connect orthogonal polynomials and combinatorics, 
we shall follow methods of set partition statistics and their 
card arrangements explained in Section \ref{sec4}.
We shall provide a combinatorial moment formula of the $s$-weighted (deformed) 
free Poisson random variables in Section \ref{sec5}. The combinatorial approach provides 
very nice combinatorial interpretations to deformation parameters, $s$ and $t$ in 
the formula. In Section \ref{sec6}, we shall investigate the limit case of $t \to 0$ with $s = 1$,
which can be regarded as the conditionally free Poissson distribution with respect 
to the reference measure of the semicircle law.  
This is based on the non-crossing partitions restricted to the case that only singletons and 
pairs are allowed to be inner.
One can see that the deformation treated in this paper interpolates
free and boolean Poisson random variables, distributions and their moments.
Moreover, one can also obtain the conditionally free Poisson distribution 
with the reference measure of the semicircle law by taking limit of 
the parameter.

%
\section{Preliminaries}\label{sec2}
\subsection{Weighted (Deformed) Free Fock Space}
Let $\hilsp$ be a real Hilbert space equipped with the inner 
product $\langle \, \cdot \,  | \, \cdot \, \rangle$, and 
$\varOmega$ be a distinguished unit vector, the so-called vacuum vector.
We denote by ${\mathcal F}_{\text{fin}} (\hilsp)$ the set of all the finite 
linear combinations of the elementary vectors 
$\xi_1 \otimes$ $\cdots$ $\otimes \xi_n$ $\in \hilsp^{\otimes n}$ \; 
$(n = 0, 1, 2, \ldots )$, where 
$\hilsp^{\otimes 0} = {\mathbb C} \varOmega$ as convention.

Let us now recall the minimum about 
the $s$-weighted free Fock space,
which is a special case of the weighted $q$-deformed Fock space 
with $q=0$ and the weight sequences $\tau_n = s^{n-1}$ $(n \ge 1)$
in \cite{AY20}\cite{BY06}.  

For $0 < s \leq 1$,  we introduce the new inner product 
$\left( \, \cdot  \, | \, \cdot \,  \right)_{s}$ 
on ${\mathcal F}_{\text{fin}} (\hilsp)$ by 
$$
  \left( \xi_1 \otimes \cdots \otimes \xi_n \,  | \,  
          \eta_1 \otimes \cdots \otimes \eta_m  \right)_{s}
  =  \delta_{m, n} \, s^{\frac{n(n - 1)}{2}} \! \! 
      \prod_{i=1}^n
       \langle \xi_i \, |  \, \eta_i \rangle,
$$
It is easy to see the positivity of the inner product
 $\big( \, \cdot \,  | \, \cdot \, \big)_{s}$.

\begin{definition}
{\it The $s$-weighted free Fock space} (simply, called {\it the $s$-free Fock space}) 
denoted by ${\mathcal F}_{s} (\hilsp)$ is
the completion of ${\mathcal F}_{\text{fin}} (\hilsp)$ with respect 
to the inner product $\left( \, \cdot \,  | \, \cdot \, \right)_{s}$. 
It is easy to see that ${\mathcal F}_{1} (\hilsp)$ is nothing but 
the free (full) Fock space  (See \cite{VDN92}, for instance).
\end{definition}

\begin{definition}
For $0 < s \leq 1$  and $\xi \in \hilsp$, 
{\it the $s$-free creation operator $a_{s}^{\dagger} (\xi)$} 
is defined by the canonical left creation,
\begin{align}\label{eq:a^*}
   &  a_{s}^{\dagger} (\xi) \, \varOmega = \xi,  \notag \\
   & a_{s}^{\dagger} (\xi) \, (\xi_1 \otimes \cdots \otimes \xi_n) 
                    = \xi \otimes \xi_1 \otimes \cdots \otimes \xi_n, 
                  \quad n \ge 1.
\end{align}
{\it The $s$-free annihilation operator $a_{s}(\xi)$}
 is defined by the adjoint operator of $a_{s}^{\dagger}(\xi)$ 
with respect to the inner product 
$\big(\, \cdot \,  | \, \cdot \, \big)_{s}$, 
that is, $a_{s}(\xi) = \Big( a_{s}^{\dagger} (\xi) \Big)^{*}$.
\end{definition}

The action of the $s$-annihilation operator on the elementary vectors 
is a direct consequence of the above definition.
\begin{proposition}
For $0 < s \leq 1$ and $\xi \in \hilsp$,  
the $s$-annihilation operator $a_{s}(\xi)$ acts on the 
  elementary vectors as follows: 
\begin{align}\label{eq:a}
  &  a_{s}(\xi) \, \varOmega = 0, \quad 
    a_{s}(\xi) \, \xi_1 = \langle \xi | \xi_{1} \rangle \, \varOmega , \notag \\ 
  & a_{s}(\xi) \, (\xi_1 \otimes \cdots \otimes \xi_n) 
      = \displaystyle{ s^{n-1} \langle \xi | \xi_1 \rangle \, 
                     \xi_2 \otimes \cdots \otimes \xi_n },
      \quad n \ge 2.
\end{align}
\end{proposition}

Moreover, let us recall other special operators on 
${\mathcal F}_{s}({\mathcal H})$.
\begin{definition}
\begin{itemize}
\item[(1)]
For $0 < s \leq 1$,  
{\it the scalar operator} $k_s$ is defined by 
 \begin{align}\label{eq:k_s}
   k_s \, & \varOmega = \varOmega, \notag \\
   k_s \, & (\xi_1 \otimes \xi_2 \otimes \cdots \otimes \xi_n) 
        =  s^n  \xi_1 \otimes \xi_2 \otimes \cdots \otimes \xi_n,
       \quad n \ge 1.
 \end{align}
\item[(2)]
For $0< t\leq 1$, 
{\it the intermediate operator} $m_t$ is defined by 
 \begin{align}\label{eq:m_t}
   m_t \, & \varOmega = 0, \notag \\
   m_t \, & (\xi_1 \otimes \xi_2 \otimes \cdots \otimes \xi_n) 
        =  t^{n-1}\xi_1 \otimes \xi_2 \otimes \cdots \otimes \xi_n, 
       \quad n \ge 1.
\end{align}

\end{itemize}
\end{definition}

\begin{remark}
We note that the operators $k_s$ for $s\in (0,1)$ and $m_t$ for $t \in (0,1)$
can be interpreted as a deformation of the identity operator $I$. 
\end{remark}

\begin{proposition}\label{prop:s-ccr}
The $s$-creation and the $s$-annihilation operators satisfy 
the following relation, 
\begin{equation*}
    a_{s}(\xi) \, a_{s}^{\dagger} (\eta) 
     = \langle \xi | \eta \rangle \, k_s,
       \qquad \xi, \eta \in \hilsp. 
\end{equation*}
\end{proposition}

A {\it noncommutative} (or {\it quantum}) probability space is a unital 
(possibly noncommutative) algebra $\mathcal{A}$ together with a linear 
functional 
$\phi : \mathcal{A} \to \mathbb{C}$, such that $\phi(1) = 1$.
If $\mathcal{A}$ is a $C^*$-algebra and $\phi$ is a state, then 
$\big( \mathcal{A}, \phi \big)$ is 
called a {\it $C^*$-probability space}.
An operator in $\mathcal{A}$ is regarded as a {\it noncommutative 
random variable} and the 
{\it distribution of $x \in \mathcal{A}$} with respect to $\phi$ is 
determined by the 
linear functional $\mu$ on $\mathbb{C}[X]$ (the polynomials in one variable) 
by 
$$
    \mu : \mathbb{C}[X] \ni P \longmapsto \phi \big( P(X) \big) \in \mathbb{C}.
$$
Considered in the $C^*$-probability context, 
the distribution $\mu$ of a 
self-adjoint operator $x \in \mathcal{A}$ can be extended to, and 
identified with the (compactly supported) probability distribution 
$\mu$ on $\mathbb{R}$ by 
$$
  \phi \big( P(X) \big) = \int_{\mathbb{R}} P(t) \, \text{d} \mu(t), 
  \quad P \in {\mathbb{C}}[X].
$$

Let us consider the vacuum state $\varphi$ for bounded operators on 
the $s$-Fock space ${\mathcal F}_{s} ({\mathscr H})$ as 
\begin{equation*}
   \varphi (b) = \left( b \, \varOmega \, | \, \varOmega \right)_{s},
   \qquad b \in \mathcal{B} \big( {\mathcal F}_{s} ({\mathscr H}) \big),
\end{equation*}
which is called {\it the vacuum expectation of $b$}.
One can employ 
$\big( \mathcal{B} \big( {\mathcal F}_{s} ({\mathscr H}) \big), 
     \varphi \big)$ 
as the noncommutative probability space, on which the model of the 
$s$-free Poisson random variable (of parameters $\lambda$ and $t$) will be discussed.

%
\section{Weighted (Deformed) Free Poisson Random Variables}\label{sec3}
From now on, let us treat the $s$-free Fock space of one-mode case 
with the unit base vector $\xi\in \hilsp$, $\|\xi\|_{\hilsp}=1$.
The $s$-creation $a_{s}^{\dagger} (\xi)$ and 
the $s$-annihilation $a_{s}(\xi)$ operators 
are simply denoted by 
$a_{s}^{\dagger}$ and $a_{s}$, respectively.
In case of one-mode, the operators $a_{s}^{\dagger}$, $a_{s}$, $k_s$, and $m_t$ 
act on the elementary vectors as follows, which can be obtained 
immediately from definitions in Section \ref{sec2}.

\begin{lemma}
For $s, t\in (0,1]$ and $\xi\in\hilsp$ with $\|\xi\|_{\hilsp}=1$, 
 \begin{align*}
     a_{s}^{\dagger} \, & \xi^{\otimes m} = \xi^{\otimes (m+1)}, \; \; m \ge 0,  \quad 
   & a_{s}           \, & \xi^{\otimes m} = 
    \begin{cases}
          s^{m-1}  \, \xi^{\otimes (m-1)}, \; & m \ge 1, \\ 
          0,                                       & m = 0,
    \end{cases} \\
     k_s \, & \xi^{\otimes m}  = s^m \, \xi^{\otimes m}, \; \; m \ge 0, \quad 
   & m_t \, & \xi^{\otimes m}  = 
    \begin{cases}
           \,t^{m-1} \xi^{\otimes m}, \; & m \ge 1, \\
          0,                           & m = 0.
    \end{cases}
 \end{align*}
\end{lemma}

By direct computations, one can see that 
the following commutation relations hold:
\begin{proposition}
For $s,t\in(0,1]$ and $\xi\in\hilsp$ with $\|\xi\|_{\hilsp}=1$,  the following equality holds:
\begin{itemize}
\item[$(1)$]
\begin{equation*}
	(a_s a_s^{\dagger})\xi^{\otimes m}
	=s(a_s^{\dagger}a_s)\xi^{\otimes m}=k_s\xi^{\otimes m}, 
	\quad m\geq 1.
\end{equation*}

\item[$(2)$]
\begin{equation*}
\begin{cases}
	\ (k_sa_{s}^{\dag})\xi^{\otimes m} = s\left(a_{s}^{\dag}k_s\right)\xi^{\otimes m}, 
	\quad m\geq 0\\
	\ s\left(k_sa_{s}\right)\xi^{\otimes m} =(a_{s}k_s)\xi^{\otimes m}, 
	\quad m\geq 1.
\end{cases}
\end{equation*}

\item[$(3)$]
\begin{equation*}
\begin{cases}
	\ (m_ta_{s}^{\dag})\xi^{\otimes m} = t\left(a_{s}^{\dag}m_t\right)\xi^{\otimes m},
	\quad m\geq 0\\
	\ t\left(m_ta_{s}\right)\xi^{\otimes m} =(a_{s}m_t)\xi^{\otimes m},
	\quad m\geq 1.
\end{cases}
\end{equation*}

\item[$(4)$]
\begin{equation*}
	(k_s m_t)\xi^{\otimes m}=(m_tk_s)\xi^{\otimes m},
	\quad m\geq 1.
\end{equation*}
\end{itemize}
\end{proposition}


\begin{definition}
For $\lambda > 0$ and $s,t\in (0,1]$,
consider a bounded self-adjoint operator 
${\bm P}^{s}_{t,\lambda}$ defined by 
\begin{equation}
  {\bm P}^{s}_{t, \lambda} 
	=  m_t +
	\sqrt{\lambda}\left(a_{s}^{\dagger} + a_{s}\right) 
	+  \lambda k_s,
\end{equation}
 on the $s$-free Fock space of one-mode.
The probability distribution of ${\bm P}^{s}_{t, \lambda}$ with 
respect to the vacuum expectation is called 
{\it the $s$-free Poisson distribution of parameters $\lambda$ and $t$}
denoted by $\Pi^{s}_{t, \lambda}$ in this paper.
\end{definition}
These operators can realize our desired model of the $s$-free Poisson 
random variables on a noncommutative probability space
$\big( \mathcal{B} \big( {\mathcal F}_{s} ({\mathscr H}) \big), 
\varphi \big)$.   
There are combinatorial meanings behind ${\bm P}^{s}_{t, \lambda}$.
It will be explained later in Section \ref{sec4} and Section \ref{sec5}.

\begin{theorem}\label{thm:s-poisson}
Suppose that $\lambda>0$ and $s, t\in (0,1]$. 
The distribution  $\Pi^{s}_{t,\lambda}$ is the orthogonalizing probability 
distribution for the sequence of orthogonal polynomials 
$\{ C_{t, n}^{s} (\lambda ; x) \}$ determined by the following recurrence relation: 
 \begin{align}\label{eq:OP-C}
    & C_{t, 0}^{s} (\lambda ; x) = 1, \quad C_{t, 1}^{s} (\lambda ; x) = x - \lambda,  \notag\\
    & C_{t, n+1}^{s} (\lambda ; x) = 
        \left( x - (\lambda  s^n +t^{n-1}) \right) C_{t, n}^{s} (\lambda ; x) 
          - \lambda s^{n-1}  \, C_{t, n-1}^{s} (\lambda ; x), \quad n \ge 1. 
 \end{align}
\end{theorem}

\begin{proof}
We simply denote the operator ${\bm P}^{s}_{t, \lambda}$ by ${\bm P}$
and $C_{t,n}^{s} (\lambda ; x)$ is abbreviated as $C_{n} (x)$, then 
it suffices to show that for $\xi \in \hilsp$
\begin{equation*}
  C_{n}({\bm P}) \varOmega = \sqrt{ \lambda^{n}} \, \xi^{\otimes n}, 
  \quad n \ge 0, 
\end{equation*}
where $\xi^{\otimes 0} = \varOmega$ because we know 
\begin{equation*}
   \varphi \big({\bm P}^n \big) = \big( {\bm P}^n \, \varOmega \, | \, \varOmega \big)_{s}
  = \int_{\mathbb{R}} x^n \, \text{d} \Pi^{s}_{t,\lambda}(x).
\end{equation*}
We shall show this by induction on $n$. It is clear that 
\begin{equation*}
  C_{0} ({\bm P}) \, \varOmega = {\mathbf 1} \, \varOmega = \varOmega,           \quad 
  C_{1} ({\bm P}) \, \varOmega = {\bm P} \, \varOmega - \lambda {\mathbf 1} \, \varOmega 
                       = \big( \sqrt{\lambda} \,  \xi + \varOmega \big) 
                                        - \lambda {\mathbf 1} \, \varOmega
                       = \sqrt{\lambda} \,  \xi.
\end{equation*}
If $n \ge 2$, we assume 
$C_{k} ({\bm P}) \, \varOmega = \sqrt{ \lambda^{k}} \, \xi^{\otimes k}$ 
for $k \le n $. Then it follows that 
 \begin{align*}
   C_{n+1}({\bm P}) \, \varOmega 
 &=  \left(
     \left( {\bm P} - (\lambda  s^n +t^{n-1})  {\mathbf 1} \right) \, C_{n}({\bm P}) 
      - \lambda  s^{n-1}  \, C_{n-1}({\bm P}) 
     \right) \,  \varOmega \\
 &=  {\bm P} \, \sqrt{ \lambda^{n}} \, \xi^{\otimes n}
    - (\lambda  s^n +t^{n-1}) \, \sqrt{ \lambda^{n}} \xi^{\otimes n} 
    - \lambda  s^{n-1} \, \sqrt{ \lambda^{n-1}} \xi^{\otimes (n-1)} \\ 
 &=  \left(m_t  + \sqrt{\lambda} \, a_{s} 
          + \sqrt{\lambda} \, a_{s}^\dagger +  \lambda k_s \right) 
                               \sqrt{ \lambda^{n}} \xi^{\otimes n} \\
 &   \qquad \qquad 
    - s^{n}   \sqrt{ \lambda^{n+2}} \xi^{\otimes n} 
    - t^{n-1} \sqrt{ \lambda^{n}  } \xi^{\otimes n} 
    - s^{n-1}1 \sqrt{ \lambda^{n+1}} \xi^{\otimes (n-1)} \\
 &= t^{n-1} \sqrt{ \lambda^{n}} \xi^{\otimes n}
    + s^{n-1} \sqrt{ \lambda^{n+1}} \xi^{\otimes n-1}
    + \sqrt{ \lambda^{n+1}} \xi^{\otimes n+1}
    + s^{n} \sqrt{ \lambda^{n+2}} \xi^{\otimes n} \\
 &   \qquad \qquad 
    - s^{n} \sqrt{ \lambda^{n+2}} \xi^{\otimes n} 
    - t^{n-1}\sqrt{ \lambda^{n}  } \xi^{\otimes n} 
    - s^{n-1} \sqrt{ \lambda^{n+1}} \xi^{\otimes (n-1)} \\
 &=   \sqrt{ \lambda^{n+1}} \xi^{\otimes n+1}.
 \end{align*}
Since $\big\{ C_n ({\bm P}) \big\}_{n \ge 0}$ are self-adjoint operators, we have 
 \begin{align*}
    \big( C_n ({\bm P}) \, C_m ({\bm P}) \, \varOmega \, | \, \varOmega \big)_{s}
& = \big( C_m ({\bm P}) \, \varOmega \, | \, C_n ({\bm P}) \, \varOmega \big)_{s} \\
& = \big( \sqrt{ \lambda^{m}} \xi^{\otimes m} \, | \, 
      \sqrt{ \lambda^{n}} \xi^{\otimes n} \big)_{s}\\
& = 0 \; \mbox{ if } \; m \ne n, 
 \end{align*}
which implies 
\begin{equation*}
  \int_{\mathbb{R}}  C_{t, n}^{s} (\lambda; x) \, C_{t,m}^{s} (\lambda; x) 
                 \, \text{d} \Pi^{s}_{t, \lambda}(x) = 0 \; \mbox{ if } \; m \ne n.
\end{equation*}
\end{proof}

\begin{remark}
One can obtain orthogonal polynomials 
of the free Poisson \cite{VDN92}  if  $s=t=1$ 
and of the boolean Poisson \cite{SW97} if $s\to 0$ and $t\to 0$ in \eqref{eq:OP-C} . 
If $s=1$ and $t\to 0$, one can yield a very interesting example, 
the conditionally free Poisson distribution.
See Section \ref{sec:s=1, t= 0}.
\end{remark}
\section{Set Partition Statistics}\label{sec4}

In our moment formula, the set partitions will be employed as combinatorial objects. 
Here we shall recall the definition of set partitions and introduce some partition 
statistics for later use.

For the set $[n] = \{1, 2, \ldots , n\}$, a {\it partition of $[n]$} is a collection
$\pi = \{ B_1, B_2, \ldots , B_k \}$ of non-empty disjoint subsets of $[n]$, 
which are called blocks and whose union is $[n]$.   
For a block $B$, we denote by $|B|$ the size of the block $B$, that is, 
the number of the elements in the block $B$. 
A block $B$ is called {\it singleton} if $|B| = 1$.
The set of all partitions of $[n]$ will be denoted by ${\mathcal P}(n)$.
The partition $\pi\in {\mathcal P}(n)$ is said to be {\it crossing} if there exist 
two blocks $B_i \ne B_j$ in $\pi$ and elements $b_1, b_2 \in B_i$, 
$c_1, c_2 \in B_j$ such that $b_1 < c_1 < b_2 < c_2$.  
A partition is called {\it non-crossing} if it is not crossing. 
We denote  by ${\mathcal NC}(n)$ the set of all non-crossing partitions of 
the set $[n]$.  
One can consult, for example,  \cite{BB15}\cite{HO07}\cite{N96}\cite{S94}\cite{SU91}\cite{Sp94}
for non-crossing partitions in detail.


\subsection{Total Depth of the Blocks by the Last and Intermediate Elements} 
For our combinatorial formula, we shall introduce the following partition 
statistics related to the last (maximum) and intermediate (neither first nor last) elements of the blocks.

For a block $C$ of the partition $\pi \in \mathcal{NC}(n)$, we consider 
the first (minimum) element $f_C$ and the last (maximum) element $\ell_C$ in 
the block $C$. In case of singleton it means $f_C = \ell_C$.
For an element $a \in [n]$, we say that the block $C$ {\it covers} $a$ 
if $a$ does not belong to the block $C$, but $a$ is included 
in the interval $[f_C, \ell_C]$.   The intermediate element $i_C$ in the block $C$ is defined to be 
neither $f_C$ nor $\ell_C$ for  $C$ with $|C|\geq 3$.  
Let us set $\text{dp}(a) = \# \{ C\in\pi \ | \ \text{$C$ covers $a$}\}$.
We note that ${\rm dp}(f_C)={\rm dp}(\ell_C)$ holds
for a block $C$ of $|C|\geq 2$ and ${\rm dp}(\ell_C)={\rm dp}(i_C)$ does 
for that of $|C|\geq 3$.

\begin{definition}
Let $B$ be a block of a non-crossing partition $\pi \in {\mathcal NC}(n)$.
\begin{itemize}
\item[(1)]
Let $\text{dp}(\ell_B)$ count the block covering 
$\ell_B$, which is called 
{\it the depth of the block $B$ by the last element}.
For $\pi\in {\mathcal NC}(n)$, the statistics $\text{td}_1(\pi)$ are 
{\it the total depth of the blocks by the last elements}
defined as 
\begin{equation*}
   \text{td}_1(\pi) = \sum_{B \in \pi, |B|\geq 1} \text{dp}(\ell_B).
\end{equation*}

\item[(2)]
Let $\text{dp}(i_B)$ count the block covering 
$i_B$,  which is called 
{\it the depth of the block $B$ by the intermediate element}.
For $\pi\in {\mathcal NC}(n)$, the statistics $\text{td}_2(\pi)$ are 
{\it the total depth of the blocks by the intermediate elements}
defined as
\begin{equation*}
   \text{td}_2(\pi) = \sum_{B \in \pi, |B|\geq 3} (|B|-2)\text{dp}(i_B). 
\end{equation*}
\end{itemize}
\end{definition}

\section{Combinatorial Moment Formula of 
the Weighted (Deformed)  Free Poisson Random Variable and Disribution}\label{sec5}
We are going to investigate the $n$-th moments of the $s$-free Poisson distribution, 
$\Pi^{s}_{t, \lambda}$. 
Namely, we evaluate the vacuum expectation of the $n$-th power of 
the $s$-free Poisson random variable (of parameters $\lambda$ and $t$), 
\begin{equation*}
 \varphi \left( \left( {\bm P}^{s}_{t, \lambda} \right)^n \right)
  = \left( \left( m_t + \sqrt{\lambda} a_{s}
           + \sqrt{\lambda} a_{s}^\dagger
               + \lambda k_s \right)^n \, \varOmega \, \Big| \, \varOmega \right)_s.
\end{equation*}
We expand 
$\left( m_t + \sqrt{\lambda} a_{s} + \sqrt{\lambda} a_{s}^\dagger
               + \lambda k_s \right)^n$ 
and evaluate the vacuum expectation in a term wise. 
In the expansion, however, we shall treat all the operators $( m_t )$, 
$( \sqrt{\lambda} \, a_{s})$, $( \sqrt{\lambda} \, a_{s}^\dagger)$, 
and $( \lambda \, k_s )$ to be noncommutative. 
A product of operators $( m_t )$, $( \sqrt{\lambda} a_s)$, 
$( \sqrt{\lambda} a_s^{\dagger})$, and $( \lambda k_s )$ are called {\it admissible}
if it has non-zero vacuum expectation.

For a given product of length $n$ 
$$
   y = z_n z_{n-1} \cdots z_2 z_1 
$$ 
where 
$$
  z_k \in \left\{ ( m_t ), \, 
                  ( \sqrt{\lambda} a_s), \, 
                  ( \sqrt{\lambda} a_s^{\dagger} ), \, 
                  ( \lambda k_s ) \right \}  \quad (k = 1, 2, \ldots ,n), 
$$
we put the sets as 
$$
 \begin{aligned}
     C_y & = \{\, k \,| \, z_k = (\sqrt{\lambda} a_s^{\dagger}) \}, &
     A_y & = \{\, k \,| \, z_k = (\sqrt{\lambda} a_s ) \}, \\ 
     M_y & = \{\, k \,| \, z_k = ( m_t ) \} , & 
     K_y & = \{\, k \,| \, z_k = (\lambda k_s) \}.
 \end{aligned}
$$
We should note that the factors are labeled from the right.  
We shall define the level of the $k$-th factor $z_k$,  
$\ell(k) \; (1 \le k \le n)$, as  
$$
   \ell(1) = 0, \quad \; 
   \ell(k) = \sum_{j=1}^{k-1} \chi(j) \quad (k \ge 2), 
$$
where $\chi(j)$ is the step function given by 
$$
 \chi(j) = \begin{cases}
                1, \; & \mbox{ if } j \in C_y, \\
               -1, \; & \mbox{ if } j \in A_y, \\
                0, \; & \mbox{ if } j \in M_y \cup K_y.
           \end{cases}
$$
Then it can be seen by rather routine argument that the monomial $y$ 
is admissible if and only if the levels $\ell(k) \; (1 \le k \le n)$ 
satisfy the following Motzkin path conditions:  
$$
 \ell(k) \ge 0   \; \, \text{\rm for } \;  1 \le k \le n, \quad 
 \ell(k) \ge 1   \: \, \text{\rm if  } \;  k \in M_y, 
 \; \; \text{and} \; \; 
 \sum_{j=1}^n \chi(j) = 0.
$$
If the monomial $y$ is an admissible product then  
the level $\ell(k)$ reflects the fact that 
$$
  \left( z_{k-1} z_{k-2} \cdots z_1 \right) \varOmega \in 
                              {\mathbb C} \xi^{\otimes \ell(k)},
$$
where $\xi^{\otimes 0}$ means the vacuum vector $\varOmega$.

It should be aware of that the operators 
$(\sqrt{\lambda} a_s^{\dagger})$ and $(\sqrt{\lambda} a_s)$ 
make a complete parenthesization in an admissible product.  
Thus we can have the non-crossing partition $\pi(y)$ in 
${\mathcal NC} (n)$ associated with an 
admissible product $y$ of length $n$ as follows: \;  

Consider the sets $C_y, A_y, M_y$, and $K_y$ as above. 
Each element in the set $K_y$ makes a singleton. The elements in 
the sets $C_y$ and $A_y$ will be used for the first and the last 
elements of blocks, respectively. The elements in the set $M_y$ will 
be used for intermediate elements of blocks.  
Of course, it will be automatically determined by non-crossingness 
that, in which block each of elements in $M_y$ should be contained, 
because the elements in the sets $C_y$ and $A_y$ are completely 
parenthesized.

\begin{example}\label{ex:adpro}
\begin{itemize}
\item[(a)]
For the admissible product of length $6$, 
$$y_a = \underbrace{ (\sqrt{\lambda} a_s  ) }_{z_6}
        \underbrace{ (\sqrt{\lambda} a_s  ) }_{z_5}
        \underbrace{ (\lambda k_s)                }_{z_4}
        \underbrace{ (\sqrt{\lambda} a_s^{\dagger}) }_{z_3}
        \underbrace{ (m_t)                        }_{z_2}
        \underbrace{ (\sqrt{\lambda} a_s^{\dagger}) }_{z_1}, 
$$
we have $C_{y_a} = \{ 1, 3 \},$ $A_{y_a} = \{ 5, 6 \},$  
$M_{y_a} = \{ 2 \},$ and $K_{y_a} = \{ 4 \} .$  
Thus we obtain the non-crossing partition, 
$$
  \pi({y_a}) = \{ \{ 1, 2, 6 \}, \{ 3,  5\} , \{ 4 \} \}.
$$

\item[(b)]
For the admissible product of length $7$, 
$$
        y_b  =
        \underbrace{ (\sqrt{\lambda} a_s  ) }_{z_7}
        \underbrace{ (\sqrt{\lambda} a_s  ) }_{z_6}
        \underbrace{ (m_t)                        }_{z_5}
        \underbrace{ (\sqrt{\lambda} a_s  ) }_{z_4}
        \underbrace{ (\sqrt{\lambda} a_s^{\dagger}) }_{z_3}
        \underbrace{ (\sqrt{\lambda} a_s^{\dagger}) }_{z_2}
        \underbrace{ (\sqrt{\lambda} a_s^{\dagger}) }_{z_1}, 
$$
we have 
$C_{y_b} = \{ 1, 2, 3 \},$ $A_{y_b} = \{ 4, 6, 7 \},$ 
$M_{y_b} = \{ 5 \},$ and $K_{y_b} = \phi.$  
Thus we obtain the non-crossing partition, 
$$
  \pi(y_b) = \{ \{ 1, 7 \}, \{ 2, 5, 6 \} , \{ 3, 4 \} \}.
$$
\end{itemize}
\end{example}

In order to evaluate the vacuum expectation of contributors, we shall use 
the cards arrangement technique which is similar as in \cite{ER96} for juggling 
patterns. We have already applied this technique to the case of non-crossing 
in \cite{YY07} and \cite{Y20}, but we are now required to prepare the different kind of 
cards. 
The cards and weights are listed below for later use.

\subsection{Creation Cards}
The creation card $C_i$  $(i \ge 0)$  has $i$ inflow lines from the 
left and $(i+1)$ outflow lines to the right, where one new line starts from the middle 
point on the ground level. For each $j \ge 1$, the inflow line of 
the $j$-th level will 
flow out at the $(j+1)$-th level without any crossing. We give the 
weight $\sqrt{\lambda}$ to the card $C_i$. 

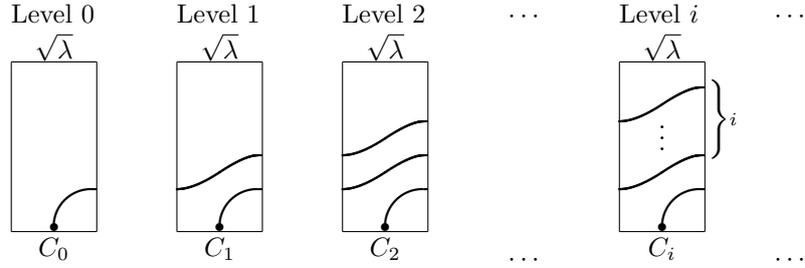
$$
 \begin{array}{llllll}
    \quad \mbox{Level } 0 & \quad \mbox{Level } 1 & 
    \quad \mbox{Level } 2 & \quad \cdots & 
    \quad \mbox{Level } i & \cdots \\
       \quad \crzer{\sqrt{\lambda}}{C_0}{} \quad 
     & \quad \crone{\sqrt{\lambda}}{C_1}{} \quad 
     & \quad \crtwo{\sqrt{\lambda}}{C_2}{} \quad 
     & \quad \cdots           \quad 
     & \quad \crnth{\sqrt{\lambda}}{C_i}{} \rsideihi{i}
     & \cdots
 \end{array} 
$$
The creation card $C_i$ represents the operation
\begin{equation*}
   \big( \sqrt{\lambda} \, a_{s}^\dagger \big) \, \xi^{\otimes i} 
     =  \sqrt{\lambda} \, \xi^{\otimes (i+1)}, \quad i \ge 0.
\end{equation*}

\subsection{Annihilation Cards}
We shall make the cards $A_i \; (i = 1,2,3, \ldots)$ for the 
$s$-annihilation operator $a_s$.  The card $A_i$ has $i$ inflow 
lines from the left and $(i-1)$ outflow lines to the right.  
On the card $A_i$, only the line of the lowest level goes 
down to the middle point on the ground level and will be annihilated.
For each $j \ge 2$, the inflow line of the $j$-th level goes 
throughout to the $(j-1)$-th level without any crossing. 
We call the card $A_i$ {\it the annihilation card of level} $i$.
We shall give the weight  $\sqrt{\lambda}s^{i-1}$ to the card $A_i$, 
where $s\in (0,1)$ counts the number of the throughout lines on the card.
It is easy to see that if $s=1$ we lose information about the depth 
of the middle point on the ground of each card.
$$
 \begin{array}{llllll}
    \quad \mbox{Level } 1 & \quad \mbox{Level } 2 & 
    \quad \mbox{Level } 3 & \quad \cdots   & 
    \quad \mbox{Level } i & \cdots  \\
       \quad \anone{\sqrt{\lambda}  }{A_1}{} \quad 
     & \quad \antwo{\sqrt{\lambda} s}{A_2}{} \quad 
     & \quad \anthr{\sqrt{\lambda} s^2}{A_3}{} \quad 
     & \quad \cdots             \quad 
     & \quad \annth{\sqrt{\lambda} s^{i-1} }{A_i}{} \rsideilo{i-1}
     & \cdots 
 \end{array}
$$
The annihilation card $A_i$ represents the operation
 \begin{equation*}
   \big( \sqrt{\lambda} \, a_{s}\big) \, \xi^{\otimes i} 
    =  \sqrt{\lambda} \, s^{i-1} \xi^{\otimes (i-1)}, 
         \qquad i \ge 1.
 \end{equation*}

\begin{remark}
A similar creation card has been used in
\cite{Y02},  but the definition is based on
the number of inner points of the arc in the block. 
In fact, the weight on the creation card $A_{i+1}$ in \cite{Y02}
is $\sqrt{\lambda}s^{2i}$.   
It is because the weight counts the number of  
inner points of the arc :
$$
  \crnth{\sqrt{\lambda}}{C_i}{} \rsideihi{i \quad \longleftrightarrow} 
  \qquad \quad 
  \lsideihi{i} \annth{ \sqrt{\lambda}s^{2i } }{A_{i+1}}{} 
$$
On the other hand, we put the weight to the annihilation card in terms of  
the depth of the last element of the block.
\end{remark}

\subsection{Scalar Cards}
The scalar card $K_i$  $(i \ge 0)$ has $i$ horizontally parallel lines and the short pole 
at the middle point on the ground. 
We shall give the weight $\lambda \, s^{i}$ to the card $K_i$, where the parameter $s\in (0,1)$ 
encodes the number of throughout lines on the card. 
It is also easy to see that if $s=1$ we lose information about the depth 
of the middle point on the ground of each card.

$$
 \begin{array}{llllll}
    \quad \mbox{Level } 0 & \quad \mbox{Level } 1 & 
    \quad \mbox{Level } 2 & \quad \cdots   & 
    \quad \mbox{Level } i & \cdots \\
       \quad \kozer{\lambda  }{K_0}{} \quad 
     & \quad \koone{\lambda s  }{K_1}{} \quad 
     & \quad \kotwo{\lambda s^2}{K_2}{} \quad 
     & \quad \cdots             \quad 
     & \quad \konth{\lambda s^{i}}{K_i}{} \rsideilo{i}
     & \cdots
 \end{array}
$$
The scalar card $K_i$ represents the operation 
\begin{equation*}
   (\lambda \, k_s) \, \xi^{\otimes i} = \lambda \, s^i \, \xi^{\otimes i},
	 \qquad i \ge 0.
\end{equation*}
For $s\in (0,1)$,  the operator $k_s$ can be viewed 
as a $s$-deformation of the identity operator $I=k_1$.

\subsection{Intermediate Cards: Non-Degenerate Case ($\bm{0<t<1}$)}

We consider the cards $M_i \; (i = 1,2,3, \ldots)$ for the operator 
$m_t$. The card $M_i$ has $i$ inflow lines and $i$, the same 
number of outflow lines.  Only the line of the lowest level goes 
down to the middle point on the ground and continues its flow as the 
lowest line again. The rest of inflow lines will keep their levels.  
We call the card $M_i$ {\it the intermediate card of level} $i$.
Since the middle point on the ground of the card 
is not the last element of the block,
we shall give a different weight $t^{(i-1)}$ from $s^{(i-1)}$ to the card $M_i$.
That is,  a parameter $t\in (0,1)$ encodes the number of directly throughout lines.
$$
 \begin{array}{llllll}
    \quad \mbox{Level } 1 & \quad \mbox{Level } 2 & 
    \quad \mbox{Level } 3 & \quad \cdots   & 
    \quad \mbox{Level } i & \cdots  \\
       \quad \nuone{1  }{M_1}{} \quad 
     & \quad \nutwo{t  }{M_2}{} \quad 
     & \quad \nuthr{t^2}{M_3}{} \quad 
     & \quad \cdots             \quad 
     & \quad \nunth{t^{(i-1)}}{M_i}{} \rsideihi{i-1}
     & \cdots 
 \end{array}
$$
The intermediate card $M_i$ represents the operation, 
$$
  m_t \xi^{\otimes i} = t^{(i-1)} \xi^{\otimes i} 
  \quad  (i \ge 1), 
$$
and the intermediate card of level $0$ is not available 
because $m_t \varOmega = 0$.

\subsection{Intermediate Cards:  Degenerate Case ($\bm{t=1}$)}
One can see easily that the case $t=1$ provides no weights to the intermediate cards $M_i$.
It means that the depth of the intermediate elements
is not counted.
Therefore,  
the intermediate cards $M_i$ with $t=1$ are labeled differently by 
$N_i \; (i = 1,2,3, \ldots)$ for the operator $m_1$ in this paper.   
$$
 \begin{array}{llllll}
    \quad \mbox{Level } 1 & \quad \mbox{Level } 2 & 
    \quad \mbox{Level } 3 & \quad \cdots   & 
    \quad \mbox{Level } i & \cdots  \\
       \quad \nuone{1  }{N_1}{} \quad 
     & \quad \nutwo{1 }{N_2}{} \quad 
     & \quad \nuthr{1}{N_3}{} \quad 
     & \quad \cdots             \quad 
     & \quad \nunth{1}{N_i}{} \rsideihi{i-1}
     & \cdots 
 \end{array}
$$
The intermediate card $N_i$ represents the operation, 
$$
  m_1 \xi^{\otimes i} = \xi^{\otimes i} 
  \quad  (i \ge 1), 
$$
and the intermediate card of level $0$ is not available 
because $m_1 \varOmega = 0$.

\begin{remark}
The degenerate intermediate card has 
been used in \cite{YY07}\cite{Y20}.
\end{remark}

\subsection{Rules for the Arrangement of the Cards:}\label{subsec:rules}

Each card arrangement gives the set partition 
of $[n]$, 
where the blocks of the partition could be obtained by the concatenation 
of the lines on the cards.  
In this construction, it is easy to find that the creation and the annihilation 
cards correspond to the first (minimum) and the last (maximum) elements in the 
blocks of size $\ge 2$, respectively, and also that the intermediate cards 
correspond to the intermediate elements in blocks.
Furthermore, {\it the weight of the arrangement} 
is given by the product of the 
weights of the cards used in the arrangement.

Now we will observe the relation between the weight of the arrangement 
and the set partition statistics:

\medskip
\noindent
{\bf On the parameter $\bm {\lambda}$:}
\begin{itemize}
\item[$\bullet$]
Since the sequence of the levels $\{ \ell(k) \}_{k = 1}^{n}$ satisfies the Motzkin 
path condition, thus we have 
$$
    \#\{ \text{creation cards} \}
  = \#\{ \text{annihilation cards} \}, 
$$
and the parameter $\lambda$ in the product of 
the weights of these cards indicates the number of the blocks of size $\ge 2$.
That is, 
$$
    \big( \sqrt{\lambda} \big)^{\#\{\text{creation cards} \}
                               +\#\{\text{annihilation cards} \} }
 = \lambda^{ \#\{ \text{creation cards}\} }
 = \lambda^{ \# \left\{ 
            B \, \mbox{\footnotesize $|$} \, B \in \pi, \, |B| \ge 2 
            \right\}}.
$$

One can see pairs of the creation 
card $C_i$ and the annihilation card $A_{i+1}$, on which the same number 
of throughout lines are drawn. 
\begin{equation}\label{eq:ca-samedepth}
  \crnth{\sqrt{\lambda}}{C_i}{} \rsideihi{i \quad \longleftrightarrow} 
  \qquad \quad 
  \lsideihi{i} \annth{ \sqrt{\lambda}s^{i } }{A_{i+1}}{} 
\end{equation}

\item[$\bullet$]
Of course, the parameter $\lambda$ in the product of the weights of scalar 
cards indicates the number of the singletons.
$$
   \lambda^{\#\{ \text{scalar cards} \}}
 = \lambda^{\# \left\{ 
            B \, \mbox{\footnotesize $|$} \, B \in \pi, \, |B| = 1 
            \right\}}.
$$
\end{itemize}
Thus the parameter $\lambda$ in the weight of an arrangement encodes 
the number of blocks of the partition,
$$
    \lambda^{ \# \left\{ 
            B \, \mbox{\footnotesize $|$} \, B \in \pi, \, |B| \ge 2 
            \right\}}
   \lambda^{\# \left\{ 
            B \, \mbox{\footnotesize $|$} \, B \in \pi, \, |B| = 1 
            \right\}}
  = \lambda^{\# \left\{ 
            B \, \mbox{\footnotesize $|$} \, B \in \pi \right\}}
  = \lambda^{|\pi|}.
$$

\medskip
\noindent
{\bf On the parameter $\bm{s}$:}
\begin{itemize}
\item[$\bullet$]
Each annihilation card corresponds to the last element $\ell_B$ of the block $B$.
In the weight of the annihilation card, the parameter $s$ counts the number 
of throughout lines on the card.  On the other hand it is clear that each throughout 
line corresponds to the block which covers the element $\ell_B$ because such a 
block should contain this throughout line as a part of concatenation.
Namely the parameter $s$ is used for encoding the depth of the block by the last 
element $\ell_B$.

\item[$\bullet$]
For the singleton $\{ k \}$, the last element of the block is itself.
Similarly, the parameter $s$ in the weight of the scalar card is used for 
encoding the number of throughout lines, which is nothing but the depth of 
the singleton $\{ k \}$.

\end{itemize}

Hence, in the weights of the card arrangement, the parameter $s\in (0,1)$ encodes 
the total depth of blocks by the last elements, that is, 
\begin{equation*}
	\left(
	\prod_{B \in \pi, |B|\geq 2}  s^{{\rm dp}(\ell_B)}
	\right)
	\left(
	\prod_{B \in \pi, |B|=1} s^{{\rm dp}(\ell_B)}
	\right)
   = s^{{\rm td}_1(\pi)} .
\end{equation*}

\begin{remark}\label{remark:dp-fl}
In case of non-crossing partitions, one can see
that ${\rm dp}(\ell_B)={\rm dp}(f_B)$ holds
for a block $B$ of $|B|\geq 2$.
See above figure in  \eqref{eq:ca-samedepth}  

\end{remark}

\medskip
\noindent
{\bf On the parameter $\bm{t}$:}
\begin{itemize}
\item[$\bullet$]
Intermediate card corresponds to the intermediate element $i_B$ of the block $B$.
In the weight of the intermediate card, the parameter $t$ counts the number 
of throughout lines on the card.  
On the other hand, it is clear that each throughout 
line corresponds to the block which covers the element $i_B$ because such a 
block should contain this throughout line as a part of concatenation.
Namely the parameter $t$ is used for encoding the depth of the block 
by the intermediate element $i_B$.
\end{itemize}

Hence, in the weights of the card arrangement, the parameter $t\in (0,1)$ encodes 
the total depth of blocks by the intermediate elements (which are not the last elements of the block), 
that is, 
\begin{equation*}
\prod_{B \in \pi, |B|\geq 3} t^{(|B|-2){\rm dp}(i_B)}
   = t^{{\rm td}_2 (\pi)} .
\end{equation*}

\begin{remark}\label{remark:dp-li}
In case of non-crossing partitions, one can see
that ${\rm dp}(\ell_B)={\rm dp}(i_B)$ holds
for a block $B$ of $|B|\geq 3$.
\end{remark}

Let $y = z_n \, z_{n-1} \cdots z_2 \, z_1$ 
be a contributor of length $n$, 
and let $\ell(k)$ the level of the $k$-th factor $z_k$.
Depending on the factors in an admissible product $y$
we shall arrange the cards along with the following rule:

If $k \in C_y$, that is, if $z_k = (\sqrt{\lambda} a_s^{\dagger})$ then 
we will put the creation card of level $\ell(k)$ with 
the $\sqrt{ \lambda }$-multiplicated weight at the $k$th position.
If $k \in A_y$, that is, if $z_k = (\sqrt{\lambda} a_s)$ then 
we will put the annihilation card of level $\ell(k)$ at the $k$-th 
position. 
The weights should be also multiplicated by $\sqrt{ \lambda }$. 
If $k \in M_y$, that is, if $z_k = (m_s)$ then we will use the 
intermediate card of level $\ell(k)$ with the original weight. 
If $k \in K_y$, that is, if $z_k = (\lambda k_s)$ 
then the singleton card of level $\ell(k)$ with the 
$\lambda$-multiplicated weight will be put at the $k$th position.
Then the non-crossing partition $\pi(y)$ can be obtained by connected 
lines because the arcs of $\pi(y)$ are naturally drawn on the cards in 
the arrangement.

Now we shall see that the vacuum expectation of an admissible 
product $y$ can be given by 
$$
  \varphi (y) = {\rm Wt} (\pi(y)).
$$
where ${\rm Wt}(\cdot)$ is  {\it the weight of the arrangement} 
given by the product of the 
weights of the cards used in the arrangement.

\begin{remark}
A similar weight function has been introduced
in \cite{Y02}\cite{Y03},  but the definition is based 
the number of inner points of the arc in the block.
See also \cite{S94} for the set partition statistics $rs$
on non-crossing partitions. 
On the other hand, we put weights by $s$ and $t$ in terms of  
the depth of the block by the last and intermediate elements, respectively.
\end{remark}

\begin{example}
\begin{itemize}
\item[(a)]
For the admissible product $y_a$ in Example \ref{ex:adpro}(a), 
we have the following card arrangement: 
$$
   \crzer{    \sqrt{\lambda}}{1}{(C_0)}
   \nuone{1                 }{2}{(M_1)}
   \crone{    \sqrt{\lambda}}{3}{(C_1)}
   \kotwo{\lambda s^2        }{4}{(K_2)}
   \antwo{\sqrt{\lambda}s }{5}{(A_2)}
   \anone{    \sqrt{\lambda}}{6}{(A_1)}
$$

The product of the cards is given by $ \lambda^3 s^3$. 
Hence, the corresponding non-crossing partition, 
$$
  \pi(y_a) = \{ \{ 1, 2, 6 \}, \{ 3,  5\} , \{ 4 \} \},
$$
has the weight ${\rm Wt} (\pi(y_a)) = \lambda^3 s^3$.

\item[(b1)] 
If we adopt the non-degenerate intermediate cards for 
the admissible product $y_b$ in the previous Example \ref{ex:adpro}(b), 
we have the following card arrangement: 
$$
 \crzer{    \sqrt{\lambda}}{1}{(C_0)}
 \crone{    \sqrt{\lambda}}{2}{(C_1)}
 \crtwo{    \sqrt{\lambda}}{3}{(C_2)}
 \anthr{\sqrt{\lambda}s^2}{4}{(A_3)}
 \nutwo{t^1                }{5}{(M_2)}
 \antwo{\sqrt{\lambda}s}{6}{(A_2)}
 \anone{    \sqrt{\lambda}}{7}{(A_1)}
$$

The product of the cards is given by $ \lambda^3 s^3t^1$. 
Hence, the corresponding non-crossing partition, 
$$
  \pi(y_b) = \{ \{ 1, 7 \}, \{ 2, 5, 6 \} , \{ 3, 4 \} \},
$$
has the weight ${\rm Wt} (\pi(y_b)) =  \lambda^3 s^3t^1$.

\item[(b2)] 
If we adopt the degenerate intermediate cards for 
the admissible product $y_b$ in Example  \ref{ex:adpro}(b), 
we have the 
following card arrangement: 
$$
 \crzer{    \sqrt{\lambda}}{1}{(C_0)}
 \crone{    \sqrt{\lambda}}{2}{(C_1)}
 \crtwo{    \sqrt{\lambda}}{3}{(C_2)}
 \anthr{\sqrt{\lambda}s^2 }{4}{(A_3)}
 \nutwo{1                 }{5}{(N_2)}
 \antwo{\sqrt{\lambda} s}{6}{(A_2)}
 \anone{    \sqrt{\lambda}}{7}{(A_1)}
$$

The product of the cards is given by $\lambda^3 s^3$. 
Hence, the corresponding non-crossing partition 
$$
  \pi(y_b) = \{ \{ 1, 7 \}, \{ 2, 5, 6 \} , \{ 3, 4 \} \}
$$
has the weight ${\rm Wt} (\pi(y_b)) = \lambda^3 s^3$.
\end{itemize}
\end{example}

Conversely, given a non-crossing partition 
$\pi \in {\mathcal NC} (n)$, 
we can make the admissible product $y(\pi)$ of length $n$ 
by the following manner: 
If $\{ k \}$ is a singleton in the partition $\pi$, 
then we put the operator $( \lambda k_s )$ as the $k$-th factor.
If $k$ is the first (resp. last) element of blocks,
then 
we use the operator $( \sqrt{\lambda} a_s^{\dagger} )$ 
(resp. $( \sqrt{\lambda} a_s)$ ) as the $k$-th factor. 
For the rest case, that is, $k$ is an intermediate element of 
a block, we adopt the operator $( m_t )$ as the $k$-th factor 
in our product. 

  Using the cards arrangement again, it is easy to see that such a 
product $y(\pi)$ has a non-zero vacuum expectation, which can be 
evaluated as the product of the weights of the cards appeared 
in the arrangement.
Hence, there is one-to-one correspondence between admissible products 
of length $n$ and the non-crossing partitions of $n$ elements 
${\mathcal NC} (n)$. 

Now we have obtained  
\begin{equation}\label{eq:moment}
  \varphi \left(( {\bm P}^s_{t, \lambda})^n \right) 
    = \! \! 
      \sum_{\begin{subarray}{c} 
               \text{\rm admissible product}  \\
               y \; \text{\rm of length } \, n 
            \end{subarray}}
       \! \! \! \! \varphi (y)  \\
    = \sum_{\pi \in {\mathcal NC} (n)} 
       \! \! {\rm Wt} (\pi).
\end{equation}
The right hand side in \eqref{eq:moment} is nothing but the 
$n$-th moment of the $s$-free Poisson distribution of parameters $\lambda$ and $t$,
$\Pi^{s}_{t, \lambda}$.  Therefore, we have derived the following combinatorial moment formula:

\begin{theorem}\label{thm:AY-moment}
Suppose $\lambda>0$ and $s,t\in (0,1]$.
The $n$-th moment of the $s$-free Poisson distribution $\Pi^{s}_{t, \lambda}$
of parameters $\lambda$ and $t$
is given by 
\begin{equation*}
  \varphi \Big( \big({\bm P}^{s}_{t, \lambda} \big)^n \Big)
    = 
	\displaystyle\sum_{\pi \in {\mathcal NC} (n)} 
       \! \!\lambda^{|\pi|} \, s^{{\rm td}_1(\pi)}t^{{\rm td}_2(\pi)},
\end{equation*}
\end{theorem}


As we mentioned in Remark \ref{remark:dp-fl} and  Remark \ref{remark:dp-li}, 
for a block $B$ in a non-crossing partition, every element in the block $B$ (regardless 
of the first, the last, or the intermediate) has the same depth. 
Hence it is possible to denote such a depth simply by $\text{dp} (B)$, 
which we call {\it the depth of the block} $B$ in a non-crossing partition. 
The notion of the depth of the block in a non-crossing partition has been treated in 
literature, for instance, \cite{AB98}\cite{HO07}.


Using the notation of $\text{dp}(B)$, 
one can have the following 
equivalent combinatorial formula:
\begin{theorem}\label{thm:AY-moment2}
Suppose $\lambda>0$ and $s,t\in (0,1]$.
The $n$-th moment of the $s$-free Poisson distribution $\Pi^{s}_{t, \lambda}$ 
of parameters $\lambda$ and $t$ is given by 
\begin{equation*}
  \varphi \Big( \big({\bm P}^{s}_{t, \lambda} \big)^n \Big)
    = 
	\displaystyle\sum_{\pi \in {\mathcal NC} (n)} \lambda^{|\pi|}
	\left\{
	\prod_{B\in \pi, |B|=1,2}  s^{{\rm dp}(B)}
	\prod_{B\in \pi, |B|\geq 3}\left(s \, t^{|B|-2} \right)^{{\rm dp}(B)}
	\right\}. 
\end{equation*}
\end{theorem}


\section{Conditionally Free Poisson with Respect to the Semicircle Law}\label{sec6}
Bo{\. z}ejko, Leinert and Speicher in \cite{BLS96} introduced notion of conditionally freeness on 
a noncommutative probability space equipped with two states, which leads to conditionally free 
convolution $\boxplus_{c}$, a binary operation on pairs of compactly supported probability 
measures on $\mathbb{R}$.  In this section, we will see that the special case of our deformed free Poisson 
yields the conditionally free Poisson distribution with the reference measure of the semicircle law.

\subsection{Conditionally Free Convolution}

In the beginning of this section, 
we shall briefly recall the definition of the conditionally 
free
convolution and the corresponding 
cumulant series in \cite{BLS96}.  

Let $\mu$ be a probability measure on $\mathbb{R}$.  Its {\it Cauchy transform}
\begin{equation*}
     G_\mu (z) = \int_{\mathbb{R}} \frac{d \mu (x)}{z - x}
\end{equation*}
is defined as as analytic function from the upper half-plane $\mathbb{C}^{+}$ into the lower-half 
plane $\mathbb{C}^{-}$ under the condition 
$\displaystyle{\lim_{y \to +\infty} i y G_{\mu}(i y) = 1}$.
If $\mu$ is a compactly supported probability measure on $\mathbb{R}$, 
then $G_{\mu}$ can be expanded into a continued fraction, 
\begin{equation}
 G_\mu (z) = \cfrac{     1}{z - \alpha_1 - 
             \cfrac{\omega_1}{z - \alpha_2 - 
             \cfrac{\omega_2}{z - \alpha_3 -  
             \cfrac{\omega_3}{z - \alpha_4 -  
             \cfrac{\omega_4}{ \; \; \ddots \; \;
             }}}}}.
\end{equation}
It should be noted that sequences $\{\omega_n\}_{n=1}^{\infty}$ of nonnegative real numbers 
and $\{\alpha_n\}_{n=1}^{\infty}$ of real numbers
are called the {\it Jacobi parameters} of  the probability 
measure $\mu$ on $\mathbb R$ 
and related to the three-term recurrence relation of the orthogonal polynomials for $\mu$.  
See \cite{Ch78}\cite{HO07} in detail, for example.


For a pair $ (\mu, \nu)$ of compactly supported probability measures $\mu$ and $\nu$ on $\mathbb{R}$, 
the free and the conditionally free cumulants series $R_{\nu}$ and $R_{(\mu,\nu)}$, respectively,  
are defined as complex functions satisfying the functional equations, 
\begin{align}
     \frac{1}{G_{\nu} (z)} & = z - R_{\nu} \big( G_\nu (z) \big), \label{eq:g_nu}\\
     \frac{1}{G_{\mu} (z)} & = z - R_{(\mu,\nu)} \big( G_\nu (z) \big). \label{eq:CFC-mu} 
\end{align}
Of course, $R_{\nu}$ is nothing but Voiculescu's $R$-transform for the free additive convolution.
Then, for pairs of compactly supported probability measures $ (\mu_1, \nu_1)$ and  $ (\mu_2, \nu_2)$,
the conditionally free convolution $(\mu, \nu) =  (\mu_1, \nu_1)  \boxplus_{c} (\mu_2, \nu_2)$ 
is defined by the requirement that both the free cumulants of the measures $\nu_i$ and 
the conditionally free cumulants of the pairs $(\mu_i, \nu_i)$ for $i=1,2$ behave additively, that is,
\begin{align*}  
     R_{\nu} (z) & = R_{\nu_1}  (z) + R_{\nu_2}  (z),  \\
     R_{(\mu, \nu)} (z) & = R_{(\mu_1, \nu_1)}  (z) + R_{(\mu_2, \nu_2)}  (z).
\end{align*} 
In particular, $\nu$ is  the Voiculescu's free convolution, $\nu_1 \boxplus \nu_2$.


\subsection{The Case of $s = 1$ and $t \to 0$}\label{sec:s=1, t= 0}

If $s = 1$ and $t \to 0$ in our deformation, 
then it can be found 
by Theorem \ref{thm:s-poisson} that 
the corresponding deformed Poisson distribution $\Pi^{1}_{0,\lambda}$  is 
the orthogonalizing probability measure for  the sequence of the polynomials 
$\big\{C_{0, k}^{1}(\lambda ; x) \big\}_{k \ge 0}$
determined by the recurrence relation,  
where $C_{0, k}^{1}(\lambda ; x) $ is simply denoted by $C_k (x)$:
\begin{equation}\label{eq:OP-C10}
\left\{
 \begin{aligned}
    & C_{0} (x) = 1, \quad C_{1} (x) = x - \lambda, \\
    & C_{2} (x) = \left( x - (\lambda  + 1) \right) C_{1} (x) - \lambda  \, C_{0}(x), \\
    & C_{n+1}(x) = 
        \left( x - \lambda \right)  C_{n}(x)  
          - \lambda \, C_{n-1}(x), \quad n \ge 2.
 \end{aligned}
\right.
\end{equation}

By using the Jacobi parameters in the recurrence relation in  \eqref{eq:OP-C10}, 
the Cauchy transform of the probability measure $\mu = \Pi^{1}_{0,\lambda}$ 
can be expressed as the continued fraction:
\begin{equation}\label{eq:CF-cfreePoi1}
 G_\mu (z) = \cfrac{     1}{z -   \lambda - 
             \cfrac{ \lambda }{z - (\lambda + 1)  - 
             \cfrac{ \lambda }{z -  \lambda  -  
             \cfrac{ \lambda }{z -  \lambda  -  
             \cfrac{ \lambda }{ \; \; \ddots \; \;
             }}}}}. 
\end{equation}
Let us put 
\begin{equation*}
 H(z) = z - \lambda - 
             \cfrac{ \lambda }{z -  \lambda  -  
             \cfrac{ \lambda }{z -  \lambda  -  
             \cfrac{ \lambda }{ \; \; \ddots \; \;
             }}}.
\end{equation*}
It is our important observation that 
the  function $H(z)$ satisfies the functional equation, 
\begin{equation}\label{eq:CF-H}
     H(z)  = z - \lambda  - \lambda \left( \sfrac{1}{H(z)} \right).
\end{equation}
Therefore the continued fraction (\ref{eq:CF-cfreePoi1}) can be rewritten as 
the form,
\begin{equation}\label{eq:CF-cfreePoi2}
 G_\mu (z) = \cfrac{ 1 }{z -   \lambda - 
                       \cfrac{ \lambda }{H(z) - 1}}
                  =  \cfrac{ 1 }{
                       z - \cfrac{ \lambda }{\quad 1 - \left( \sfrac{1}{H(z)} \right) \quad }}.
\end{equation}
On the other hand,  if $\nu$ is the semicircle law of mean $\lambda$ and 
variance $\lambda$, then the free cumulant series (Voiculescu's $R$-transform) of $\nu$
 is known to be 
\begin{equation*}
   R_{\nu} (z) = \lambda + \lambda z.
\end{equation*}
Hence, the equation of $G_\nu (z)$ in \eqref{eq:g_nu}
is equivalent to  
\begin{equation*}
     \frac{1}{G_{\nu} (z)} = z - (\lambda + \lambda  G_\nu (z)).
\end{equation*}
This equation implies that the function $H(z)$ is the reciprocal of  $G_\nu (z)$ 
by \eqref{eq:CF-H}.
Thus by substituting $G_\nu (z)$  for $\dfrac{1}{H(z)}$ into the right hand side of 
\eqref{eq:CF-cfreePoi2}, 
one can obtain 
\begin{equation}\label{eq:CT-mu} 
    \frac{1}{G_\mu(z)} = z - \frac{\lambda}{1 -G_\nu (z)}.
\end{equation}
By comparing \eqref{eq:CT-mu} with \eqref{eq:CFC-mu}, 
one can see that the conditionally 
free cumulant series $R_{(\mu, \nu)}(z)$ is given by 
\begin{equation*}
   R_{(\mu, \nu)}(z) = \frac{\lambda}{1 - z} 
           = \lambda +  \lambda \, z + \lambda \, z^2 + \lambda \, z^3 + \cdots.
\end{equation*}
This feature is consistent with the characterization of the Poisson distribution, 
that is, the constant cumulants 
of all orders.  Therefore, one can claim the following characteristic:


\begin{proposition}
The probability measure $\mu = \Pi^{1}_{0,\lambda}$ can be regarded as 
the conditionally free Poisson distribution with respect to 
the reference measure $\nu$ of the semicircle law with mean $\lambda$ 
and variance $\lambda$.
\end{proposition}

\subsection{Remarks on the Moments of the Conditionally Free Poisson Distribution}

Based on the depth of blocks in non-crossing partitions, the notion of outer or inner for the blocks in 
non-crossing partition can be introduced.  See \cite{BLS96}, for instance.

\begin{definition}
The block $B$ of a non-crossing partition is called {\it outer} if $\text{dp}(B) = 0$ and 
 {\it inner} if  $\text{dp}(B) \ge 1$.
\end{definition}

If taking limit $t \to 0$ in the moment formula in Theorem \ref{thm:AY-moment2}, 
one can see that if the size of block $\ge 3$, then only the outer (of depth $0$) 
case will survive for the summation.  Namely, no block of size $\ge 3$ is allowed to 
be inner.

Based on the above, we will consider some restricted class of non-crossing partitions, 
that is, the set of non-crossing partitions  
with exactly $k$ blocks such that only singletons (blocks of size $1$) and pairs (blocks of size $2$) 
are allowed to be inner.  From now on, we denote it by $\mathcal{NC}_{1,2 : \text{inner}}(n,k)$.

Now we can state the following combinatorial moment formula for the case, $s=1$ and $t \to 0$.
\begin{proposition}
Suppose $\lambda > 0$.  
The $n$-th moment of $\mu = \Pi^{1}_{0,\lambda}$,
the conditionally free Poisson distribution with the reference measure $\nu$ 
being the semicircle law of mean $\lambda$ and variance $\lambda$, is given by 
\begin{equation}\label{eq:CMform-CFP}
   m_n(\mu) 
    =  
      \sum_{\pi \in \mathcal{NC}(n)}
      \Bigg(
             \prod_{B \in \pi, |B|=1,2} \! \!\! \! \! \! \lambda  \; \; 
      \Bigg)
      \Bigg( 
                \prod_{\begin{subarray}{c} 
                             B \in \pi, \, |B| \ge 3  \\
                             B :  \text{outer} 
                              \end{subarray}}
                \! \! \! \! \lambda \; \;
       \Bigg)
    =  
      \sum_{\pi\in \mathcal{NC}(n)}
      \Bigg( 
                \prod_{\begin{subarray}{c} 
                             B \in \pi, \, |B|=1, 2\\
                             B :  \text{inner} 
                              \end{subarray}}
                \! \! \! \! \lambda \; \;
       \Bigg)
      \Bigg(
             \prod_{B \in \pi,  B : \text{outer}} \! \!\! \! \! \! \lambda  \; \; 
      \Bigg),
\end{equation}
where we let $m_n(\mu)=\varphi \left( \big({\bm P}^{1}_{0, \lambda} \big)^n \right)$.
It can be written equivalently by
\begin{equation*}
        m_n (\mu) = \sum_{k = 1}^n
 {\#} \Big( \mathcal{NC}_{1,2 : \text{inner}}(n,k) \Big)  \lambda^k.
\end{equation*}
\end{proposition}

\begin{remark}
The combinatorial formula (\ref{eq:CMform-CFP}) can be also derived from 
the conditionally free convolution in Section 3 of \cite{BLS96}.
\end{remark}

\begin{example}
Here we shall list the moments of $\mu = \Pi^{1}_{0,\lambda}$ for the first few orders below:
\begin{align*}
 m_1 (\mu) &= \lambda,\\
 m_2 (\mu) &= \lambda^2 +    \lambda,\\
 m_3 (\mu) &= \lambda^3 +  3 \lambda^2 +   \lambda,\\
 m_4 (\mu) &= \lambda^4 +  6 \lambda^3 +6  \lambda^2 +    \lambda,\\
 m_5 (\mu) &= \lambda^5 + 10 \lambda^4 +20 \lambda^3 +  9 \lambda^2 + \lambda,\\
 m_6 (\mu) &= \lambda^6 + 15 \lambda^5 +50 \lambda^4 + 44 \lambda^3 + 12 \lambda^2 + \lambda,\\
m_7(\mu)  & = \lambda^7 + 21 \lambda^6 + 105 \lambda^5 + 154 \lambda^4 + 
                     77 \lambda^3 + 15 \lambda^2 + \lambda.
\end{align*}
The sequence of the moments for $\lambda = 1$ is 
\begin{equation}\label{eq:mmt_sqc}
\{1, 2, 5, 14, 41, 123, 374, 1147, 3538, 10958, \ldots \}.
\end{equation}
This sequence can be found in \cite{Sl23}  as A3262548 (the number of non-capturing set partitions of $[n]$).
See also A054391 in \cite{Sl23}.
\end{example}


\begin{remark}
Since the Cauchy transform $G_{\mu} (z)$ for $\mu = \Pi^{1}_{0,\lambda}$ satisfies 
the quadratic equation,
\begin{equation*}
       \big( z^3 - (1 + 3 \lambda)  z^2 + 3 \lambda^2  z - \lambda^3 \big)  \big( G_{\mu} (z) \big)^2
     - \big(2 z^2 -(2 + 5 \lambda) z + 3 \lambda^2 \big) G_{\mu} (z)
    + \big( z - (2 \lambda+1) \big) = 0,
\end{equation*}
one can solve it explicitly as 
\begin{equation}\label{eq:G mmt sqc}
  G_\mu(z) = 
         \frac{2 z^2 - (2 + 5 \lambda) z + 3 \lambda^2 +
                      \lambda \sqrt{ \; (z - \lambda)^2 - 4 \lambda  \;}  \; }
                 { 2 \big (z^3 - (1 + 3 \lambda) z^2 + 3 \lambda^2 z  - \lambda^3  \big)},
\end{equation}
where the branch of the square root is, of course,  chosen so that $G_\mu(z)$ is continuous 
for $z \in \mathbb{C}^{+}$.
Applying the above formula \eqref{eq:G mmt sqc}, 
the generating function $M(z)$ of the sequence \eqref{eq:mmt_sqc} is given by 
\begin{equation*}
    M(z) = \left. 
                 \sfrac{1}{\, z \,} G_\mu \Big( \sfrac{1}{\, z \,}  \Big) 
               \right|_{\lambda = 1}
            =\frac{-3 z^2 +7 z - 2 - z \sqrt{-3 z^2-2 z+1}}
                        {2 \left(z^3 - 3 z^2 + 4 z - 1 \right)}.
\end{equation*}
\end{remark}


\begin{remark}[Boolean and Fermionic cases]
\; \vspace{-3pt}
\begin{itemize}
\item[(1)]
In case of $s \to 0$ and $t \to 0$, the non-crossing partitions in the combinatorial formula in 
Theorem \ref{thm:AY-moment2}  for $\Pi^{0}_{0,\lambda}$ are restricted to the case that 
all the blocks are of depth $0$, that is, no inner block is allowed.  Such non-crossing partitions 
are called the {\it interval partitions}  ($\mathcal{IP}$), which correspond to the boolean case.
\item[(2)]
In \cite{Or02},  the non-crossing partitions, 
where only singletons are allowed to be inner, were investigated.  
These are called the {\it almost interval partitions} 
($\mathcal{AIP}$).   Associated with the $\mathcal{AIP}$, 
the fermi convolution was introduced 
and the corresponding fermionic Poisson distribution discussed in \cite{SY00}
was derived.  It is claimed \cite{Or02} that  fermionic Poisson distribution
is different from boolean Poisson distribution \cite{SW97}.
\end{itemize}
\end{remark}


%


\begin{thebibliography}{99999999}
\bibitem
{AB98}
L. Accardi, M. Bo{\.z}ejko, 
{Interacting Fock spaces and gaussianization of probability measures}, 
{\em Infin. Dimens. Anal. Quantum Probab. Relat. Top.}, 
{\bf  1}(4) (1998), 663–670.
%
%
\bibitem
{AY20}
N. Asai and H. Yoshida,
{Deformed Gaussian operators on weighted $q$-Fock spaces},
{\em J. Stoch. Anal.},
{\bf 1}(4) (2020), Article 6.

\bibitem
{B01}
M. Bo{\.z}ejko,
{Deformed free probability of Voiculescu},
 {\em RIMS Kokyuroku},
{\bf 1227} (2001),  96-113.

\bibitem
{BB15}
M. Bo\.{z}ejko and W. Bo\.{z}ejko,
{Generalized Gaussian processes and relations 
with random matrices and positive definite functions on permutation groups}, 
{\em Infin. Dimens. Anal. Quantum Probab. Relat. Top.}, 
{\bf  18}(3) (2015), 1550020, 19 pp.

\bibitem
{BKW06}
M. Bo\.{z}ejko,  A.D. Krystek, and \L.J. Wojakowski,
{Remarks on the $r$  and $\Delta$ convolutions}, 
{\em Math. Z.}, {\bf 253}(1), (2006), 177-196.

\bibitem
{BLS96}
M. Bo\.{z}ejko, M. Leinert, and R. Speicher,
{Convolution and limit theorems for conditionally free random variables}, 
{\em Pacific J. Math.}, 
{\bf 175} (1996), 357--388.
%
\bibitem
{BY06}
M. Bo\.{z}ejko and H. Yoshida,
{Generalized $q$-deformed Gaussian random variables},
{\em Banach Center Publ.},
{\bf 73}(1) (2006), 127--140.

\bibitem
{Ch78}
T. S. Chihara,
{\em  An Introduction to Orthogonal Polynomials}, 
Gordon and Breach Science publishers,
New York, NY, 1978.

\bibitem
{ER96}
R. Ehrenborg and M. Readdy, 
{Juggling and applications to $q$-analogues}, 
{\em Discr. Math.}, {\bf 157}(1--3)  (1996), 107--125.

\bibitem
{HO07}
A. Hora and N. Obata,
{\em  Quantum Probability and Spectral Analysis
of Graphs}, 
Theoretical and Mathematical Physics,
Springer-Verlag, Heidreberg, 2007

\bibitem
{N96}
A. Nica,
{$R$-transforms of free joint distributions 
and non-crossing partitions},
{\em J. Funct. Anal.}, 
{\bf 135} (1996), 271--296.
%
%
\bibitem
{Or02}
F. Oravecz,
{Fermi convolution},
{\em Infin. Dimens. Anal. Quantum Probab. Relat. Top.},
{\bf 5} (2002), 235--242.
%
%
\bibitem
{SY00}
N. Saitoh and H. Yoshida,
{$q$-deformed Poisson random variables on $q$-Fock space},
{\em J. Math. Phys.},
{\bf 41}(8) (2000), 5767--5772.
\bibitem
{S94}
R. Simion,
{Combinatorial statistics on non-crossing partitions},
{\em  J. Combin. Theory,  Ser. A}, 
{\bf 66} (1994), 270--301.
%
%
\bibitem
{SU91}
R. Simion and D. Ullman, 
{On the structure of the lattice of non-crossing partitions}, 
{\em Discr. Math.},
{\bf 98}(3) (1991), 193--206.

%
%
\bibitem
{Sl23}
 N.J.A. Sloane, 
{\em The On-Line Encyclopedia of Integer Sequences}, 
2023,  available at  https://oeis.org/.
%
%

\bibitem
{Sp90}
R. Speicher,
{A new example of 'Independence' and 'White Noise'}, 
{\em Probab. Theory  Relat. Fields}, 
{\bf 84} (1990), 141--159.
%

\bibitem
{Sp94}
R. Speicher,
{Multiplicative functions on the lattice of non-crossing partitions 
and free convolution}, 
{\em Math. Ann.}, 
{\bf 298} (1994), 611--628.

%
\bibitem
{SW97}
R. Speicher and R. Woroudi,
Boolean convolution,
in: {\em Free Probability Theory}, 
D. Voiculescu (ed.), pp. 267--280,
Filelds Inst. Commun. {\bf 12}, 
Amer. Math. Soc.,  Privince RI,
1997.
\bibitem
{VDN92}
D. Voiculescu, K. Dykema, and A. Nica,
{\em  Free Random Variables. 
A Noncommutative Probability Approach to Free Products with 
Applications to Random Matrices, Operator Algebras and Harmonic 
Analysis on Free Groups}, 
CRM Monograph Series, vol.{\bf 1}.  Amer. Math. Soc., 
Providence, RI,  1992

\bibitem
{YY07}
F. Yano and H. Yoshida, 
{Some set partition statistics in non-crossing partitions and 
generating functions}, 
{\em Discr. Math.}, {\bf 307}(24) (2007), 3147-3160.

\bibitem
{Y02} 
H. Yoshida,
Remarks on the $s$-free convolution,
in: {\em Non-commutativity, Infinite Dimensionality,
and Probability at the Crossroads}, N. Obata et al. (eds.),
World scientific, Singapore,
2002, 412--433.

\bibitem
{Y03}
H. Yoshida,
{The weight function on non-crossing partitions for 
the $\varDelta$-convolution}, 
{\em Math. Z.}, {\bf 245}(1), (2003), 105-121.

\bibitem
{Y20}
H. Yoshida, 
{Remarks on a free analogue of the beta prime distribution},
{\em J. Theoret. Probab.}, 
{\bf 33}(3) (2020), 1363--1400.

\end{thebibliography}
\end{document}